\documentclass[11pt,twoside]{article}
\usepackage{amsfonts}
\usepackage{amsmath}
{\markboth{\sl  C. Abdelkefi and M.
Rachdi}{\sl
 Riesz potentials in Dunkl analysis}
\newtheorem{thm}{Theorem}[section]
\newtheorem{rem}{Remark}[section]
\newtheorem{ex}{Example}[section]
\newtheorem{prop}{Proposition}[section]
\newtheorem{cor}{Corollary}[section]

\numberwithin {equation}{section}
\newenvironment{pf}{\textbf{Proof.}} {\hfill $\Box$}
\begin{document}
\title{\LARGE\bf  Some properties of the Riesz potentials in Dunkl analysis}
\author{Chokri Abdelkefi and Mongi Rachdi
 \footnote{\small This work was
completed with the support of the DGRST research project LR11ES11
and the program CMCU 10G / 1503.}\\ \small  Department of
Mathematics, Preparatory
Institute of Engineer Studies of Tunis  \\ \small 1089 Monfleury Tunis, University of Tunis, Tunisia\\
  \small E-mail : chokri.abdelkefi@ipeit.rnu.tn \\  \small E-mail : rachdi.mongi@ymail.com}%
\date{}
\maketitle
\begin{abstract} In Dunkl theory on $\mathbb{R}^d$ which generalizes classical Fourier analysis, we study first the behavior at infinity of the Riesz
potential of a non compactly supported function. Second, we give for
$1 < p \leq q < +\infty$, weighted $Lp\rightarrow Lq$ boundedness of
the Riesz potentials with sufficient conditions. As application, we
prove a weighted generalized Sobolev inequality.
\end{abstract}
%% ----------------------------------------------------------------------
\section{Introduction}
\par
Dunkl operators $T_i, 1 \leq i \leq d$ introduced by C.F. Dunkl in
[9], are differential-difference operators associated with a finite
reflection group $G$, acting on some Euclidean space. These
operators attached with a positive root system $R_+$ and a non
negative multiplicity function $k$, can be considered as
perturbations of the usual partial derivatives by reflection parts.
They provide a useful framework for the study of multivariable
analytic structures which reveal certain reflection symmetries.
During the last years, these operators have gained considerable
interest in various fields of mathematics (see [1, 2, 3, 14, 18])
and also in physical applications; they are, for example, naturally
connected with certain Schr\"{o}dinger operators for
Calogero-Sutherland-type quantum many body systems (see [21]). The
Dunkl kernel $E_k$ has been introduced by C.F. Dunkl in [10]. For a
family of weight functions $w_k$ invariant under the action of $G$,
we use the Dunkl kernel and the measure $d\nu_k(x)=w_k(x)dx$ to
define the Dunkl transform $\mathcal{F}_k$, which enjoys properties
similar to those of the classical Fourier transform $\mathcal{F}$.
If the parameter $k\equiv0$ then $w_k(x)=1$ and the measure $\nu_k$
coincide with the Lebesgue measure, so that $\mathcal{F}_k$ becomes
$\mathcal{F}$ and the $T_i, 1 \leq i \leq d$ reduce to the
corresponding partial derivatives $\frac{\partial}{\partial x_i}, 1
\leq i \leq d$. Therefore Dunkl analysis can be viewed as a
generalization of classical Fourier analysis (see next section,
Remark 2.1). The classical Fourier transform behaves well with the
translation operator $f \mapsto f(.-y)$, which leaves the Lebesgue
measure on $\mathbb{R}^d$ invariant. However, the measure $w_k(x)dx$
is no longer invariant under the usual translation. One ends up with
the Dunkl translation operators $\tau_x$, $x\in \mathbb{R}^d$,
introduced by K. Trim\`eche in [20] on the space of infinitely
differentiable functions on $\mathbb{R}^d$ (see next section).

Let $0<\alpha < d $. The operator\begin{eqnarray*} I_{\alpha}f(x) =
\int_{\mathbb{R}^{d}}\frac{f(y)}{\|x-y\|^{d-\alpha}}dy,
\end{eqnarray*}
is known as the classical Riesz potential of $f$. If $\alpha=2$,
$d>2$ and $f$ is H\"{o}lder continuous, we have the Newtonian
potential and $I_{\alpha}f$ is a classical solution of the Poisson
equation, $\Delta u=-f$ in $\mathbb{R}^{d}$. The special case
$\alpha=2$ and $d=3$, gives the electric potential of a statistic
charge distribution with charge $f$. the behavior of the Riesz
potential of function at infinity was investigated in [4, 13]. It is
easy to see that if $f$ is non negative and compactly supported,
then $I_{\alpha}f(x)$ has the order $\|x\|^{\alpha-d}$ at infinity.
D.Siegel and E.Talvila [17] found necessary and sufficient
conditions on $f$ for the validity of
$I_{\alpha}f(x)=O(\|x\|^{\alpha-d})$ as $\|x\|\rightarrow +\infty $
even when $f$ is not compactly supported.

Our aim is first to extend these results to the context of Dunkl
theory where a similar operator is already defined on the Schwartz
space $S(\mathbb{R}^{d})$. For $0<\alpha<2\gamma+d$ with
$\displaystyle\gamma= \sum_{\xi \in R_+}k(\xi)$, the Riesz potential
$I_{\alpha}^{k}f$ of a function $f$ (see [5, 12, 19]), is given by
\begin{eqnarray*}
I_{\alpha}^{k}f(x)=2^{\gamma+\frac{d}{2}-\alpha}\frac{\Gamma(\gamma+\frac{d-\alpha}{2})}{\Gamma(\frac{\alpha}{2})}
\int_{\mathbb{R}^{d}} \frac{\tau_{y}f(x)}{\|y\|^{2\gamma+d-\alpha}}
d\nu_{k}(y),
\end{eqnarray*}
where $\nu_k$ is the weighted measure defined by
\begin{eqnarray*}d\nu_k(x):=w_k(x)dx\quad
\mbox{with}\;\;w_k(x) = \prod_{\xi\in R_+} |\langle
\xi,x\rangle|^{2k(\xi)}, \quad x \in \mathbb{R}^d.\end{eqnarray*}
$\langle .,.\rangle$ being the standard Euclidean scalar product on
$\mathbb{R}^d$ (see next section).

Second, we give for $1<p\leq q<+\infty$, sufficient conditions on
the decreasing rearrangement of non-negative locally integrable
weight functions $u$, $v$ on $\mathbb{R}^d$, such that the Riesz
potential $I_{\alpha}^{k}$ satisfies the weighted inequality
\begin{eqnarray*}\Big(\int_{\mathbb{R}^{d}}|I_{\alpha}^{k}f(y)|^{q}u(y)
d\nu_{k}(y)\Big)^{\frac{1}{q}} \leq
c\,\Big(\int_{\mathbb{R}^{d}}|f(x)|^{p}v(x)
d\nu_{k}(x)\Big)^{\frac{1}{p}},\end{eqnarray*} for $f\in
L^p_{k,v}(\mathbb{R}^d)$ with $L^p_{k,v}(\mathbb{R}^d)$ is the space
$L^{p}(\mathbb{R}^d, v(x) d\nu_k(x))$.

As consequence, we obtain weighted $L^p\rightarrow L^q$ boundedness
of the fractional maximal operator. Finally, we prove a weighted
generalized Sobolev inequality. These are generalizations of some results obtained in [5, 12]. \\

The contents of this paper are as follows. \\In section 2, we
collect some basic definitions and results about harmonic analysis
associated with Dunkl operators .\\
We study in section 3, the behavior at infinity of the Riesz
potential of a non compactly supported function.
\\In section 4, we give for $1 < p \leq q < +\infty$, weighted $Lp\rightarrow Lq$ boundedness
of the Riesz potentials with sufficient conditions. As consequence,
we obtain weighted inequalities for the fractional
maximal operators and we prove a weighted generalized Sobolev inequality. \\

Along this paper, we denote by $\langle .,.\rangle$, the standard
Euclidean scalar product on $\mathbb{R}^d$ and we write for $x \in
\mathbb{R}^d, \|x\| = \sqrt{\langle x,x\rangle}$. We use $c$ to
denote a suitable positive constant which is not necessarily the
same in each occurrence. Furthermore, we denote by

$\bullet\quad \mathcal{E}(\mathbb{R}^d)$ the space of infinitely
differentiable functions on $\mathbb{R}^d$.

$\bullet\quad \mathcal{S}(\mathbb{R}^d)$ the Schwartz space of
functions in $\mathcal{E}( \mathbb{R}^d)$ which are rapidly
decreasing as well as their derivatives.

$\bullet\quad \mathcal{D}(\mathbb{R}^d)$ the subspace of
$\mathcal{E}(\mathbb{R}^d)$ of compactly supported functions.
\section{Preliminaries}
$ $ In this section, we recall some results in Dunkl
theory (see[8, 9, 10, 16]) and we refer for more details to the surveys [15].\\

Let $G\subset O(\mathbb{R}^{d})$ be a finite reflection group on
$\mathbb{R}^{d}$, associated with a root system $R$. For $\alpha\in
R$, we denote by $\mathbb{H}_\alpha$ the hyperplane orthogonal to
$\alpha$. For a given
$\beta\in\mathbb{R}^d\backslash\bigcup_{\alpha\in R}
\mathbb{H}_\alpha$, we fix a positive subsystem $R_+=\{\alpha\in R:
\langle \alpha,\beta\rangle>0\}$. We denote by $k$ a nonnegative
multiplicity function defined on $R$ with the property that $k$ is
$G$-invariant. We associate with $k$ the index
$$\gamma = \sum_{\xi \in R_+} k(\xi),$$
and a weighted measure $\nu_k$ given by
\begin{eqnarray*}d\nu_k(x):=w_k(x)dx\quad
\mbox{ where }\;\;w_k(x) = \prod_{\xi\in R_+} |\langle
\xi,x\rangle|^{2k(\xi)}, \quad x \in \mathbb{R}^d,\end{eqnarray*}
%%%%%%%%%%%%%%%%%%%%%%%%%%%%%%%%%%%%%%%%%%%%%%%%%%%%%%%%%%%%%%%%%%%%%%%%%%%%%%%%%%%%%%%%%%%%%%%%%%%%%%%%%%%%%%%%%%PAGE 3%%%%%%%%%%%%%%%%%%%%%%%%%%%%%%%%%%%%%%

Further, we introduce the Mehta-type constant $c_k$ by
$$c_k = \left(\int_{\mathbb{R}^d} e^{- \frac{\|x\|^2}{2}}
w_k (x)dx\right)^{-1}.$$

For every $1 \leq p \leq + \infty$, we denote by
$L^p_k(\mathbb{R}^d)$,
 the spaces $L^{p}(\mathbb{R}^d,
d\nu_k(x)),$
 and we use $\| \;\|_{p,k}$ as a shorthand for $\|\ \;\|_{L^p_k( \mathbb{R}^d)}$. \\
By using the homogeneity of degree $2\gamma$ of $w_k$, for a radial
function $f$ in $L^1_k ( \mathbb{R}^d)$, there exists a function $F$
on $[0, + \infty)$ such that $f(x) = F(\|x\|)$, for all $x \in
\mathbb{R}^d$. The function $F$ is integrable with respect to the
measure $r^{2\gamma+d-1}dr$ on $[0, + \infty)$ and we have
 \begin{eqnarray} \int_{\mathbb{R}^d}  f(x)\,d\nu_k(x)&=&\int^{+\infty}_0
\Big( \int_{S^{d-1}}f(ry)w_k(ry)d\sigma(y)\Big)r^{d-1}dr\nonumber\\
&=&
 \int^{+\infty}_0
\Big( \int_{S^{d-1}}w_k(ry)d\sigma(y)\Big)
F(r)r^{d-1}dr\nonumber\\&= & d_k\int^{+ \infty}_0 F(r)
r^{2\gamma+d-1}dr,
\end{eqnarray}
  where $S^{d-1}$
is the unit sphere on $\mathbb{R}^d$ with the normalized surface
measure $d\sigma$  and \begin{eqnarray*}d_k=\int_{S^{d-1}}w_k
(x)d\sigma(x) = \frac{c^{-1}_k}{2^{\gamma +\frac{d}{2} -1}
\Gamma(\gamma + \frac{d}{2})}\;.  \end{eqnarray*}

The Dunkl operators $T_j\,,\ \ 1\leq j\leq d\,$, on $\mathbb{R}^d$
associated with the reflection group $G$ and the multiplicity
function $k$ are the first-order differential- difference operators
given by
$$T_jf(x)=\frac{\partial f}{\partial x_j}(x)+\sum_{\alpha\in R_+}k(\alpha)
\alpha_j\,\frac{f(x)-f(\rho_\alpha(x))}{\langle\alpha,x\rangle}\,,\quad
f\in\mathcal{E}(\mathbb{R}^d)\,,\quad x\in\mathbb{R}^d\,,$$ where
$\rho_\alpha$ is the reflection on the hyperplane
$\mathbb{H}_\alpha$ and $\alpha_j=\langle\alpha,e_j\rangle,$
$(e_1,\ldots,e_d)$ being the canonical basis of
$\mathbb{R}^d$.
%%%%%%%%%%%%%%%%%%%%%%%%%%%%%%%%%%%%%%%%%%%%%%%%%%%%%%%%%%%%%%%%%%%%%%%%%%%%%%%%%%%%%%%%%%%%%%%%%%%%%%%%%%%%%%%%%%%%PAGE4%%%%%%%%%%%%%%%%%%%%%%%%%%%%%%%%%%%%
\begin{rem}In the case $k\equiv0$, the weighted function $w_k\equiv1$ and the measure $\nu_k$ associated to the
Dunkl operators coincide with the Lebesgue measure. The $T_j$ reduce
to the corresponding partial derivatives. Therefore Dunkl analysis
can be viewed as a generalization of classical Fourier
analysis.\end{rem}

For $y \in \mathbb{C}^d$, the system
$$\left\{\begin{array}{lll}T_ju(x,y)&=&y_j\,u(x,y),\qquad1\leq j\leq d\,,\\  &&\\
u(0,y)&=&1\,.\end{array}\right.$$ admits a unique analytic solution
on $\mathbb{R}^d$, denoted by $E_k(.,y)$ and called the Dunkl
kernel. This kernel has a unique holomorphic extension to
$\mathbb{C}^d \times \mathbb{C}^d $.

The Dunkl transform $\mathcal{F}_k$ is defined for $f \in
\mathcal{D}( \mathbb{R}^d)$ by
$$\mathcal{F}_k(f)(x) =c_k\int_{\mathbb{R}^d}f(y) E_k(-ix, y)d\nu_k(y),\quad
x \in \mathbb{R}^d.$$  We list some known properties of the Kernel
$E_k$ and the Dunkl transform:
\begin{itemize}
\item[i)] For all $\lambda\in
\mathbb{C}$ and $z, z'\in \mathbb{C}^d,$ we have\\ $E_k(z,z') =
E_k(z',z)$; \;$E_k(\lambda z,z') = E_k(z,\lambda z').$
 \\For $x, y \in
\mathbb{R}^d,$  $|E_k(x,iy)| \leq 1$.
\item[ii)] The Dunkl transform of a function $f
\in L^1_k( \mathbb{R}^d)$ has the following basic property
\begin{eqnarray*}\| \mathcal{F}_k(f)\|_{\infty,k} \leq
 \|f\|_{ 1,k}\;. \end{eqnarray*}
\item[iii)] The Dunkl transform is an automorphism on the Schwartz space $\mathcal{S}(\mathbb{R}^d)$.
\item[iv)] When both $f$ and $\mathcal{F}_k(f)$ are in $L^1_k( \mathbb{R}^d)$,
 we have the inversion formula \begin{eqnarray*} f(x) = \int_{\mathbb{R}^d}\mathcal{F}_k(f)(y) E_k( ix, y)d\nu_k(y),\quad
x \in \mathbb{R}^d.\end{eqnarray*}
\item[v)] (Plancherel's theorem) The Dunkl transform on $\mathcal{S}(\mathbb{R}^d)$
 extends uniquely to an isometric automorphism on
$L^2_k(\mathbb{R}^d)$.
\item[vi)] For $f\in\mathcal{S}(\mathbb{R}^d)$ and $1\leq j \leq d $, we have
\begin{eqnarray}\mathcal{F}_{k}(T_{j}f)(\xi)=i\xi_{j}\mathcal{F}_{k}(f)(\xi), \;\;\xi\in\mathbb{R}^{d}.\end{eqnarray}
\end{itemize}
The Dunkl translation operator $\tau_x$, $x\in\mathbb{R}^d$, was
introduced in [20] on $\mathcal{E}( \mathbb{R}^d)\;.$ For $f\in
\mathcal{S}( \mathbb{R}^d)$ and $x, y\in\mathbb{R}^d$, we have
\begin{eqnarray*}\mathcal{F}_k(\tau_x(f))(y)=E_k(i x, y)\mathcal{F}_k(f)(y).\end{eqnarray*}  The space
\begin{eqnarray*}A_k(\mathbb{R}^d)=\{f\in L^1_k(\mathbb{R}^d):\;\mathcal{F}_{k}(f)\in L^1_k(\mathbb{R}^d)\}. \end{eqnarray*}
plays a particular role in the analysis of this generalized
translation (see [16, 18, 20]). Observe that $A_k(\mathbb{R}^d)$ is
contained in the intersection of $L^1_k(\mathbb{R}^d)$ and
$L^{\infty}$ and hence is a subspace of $L^2_k(\mathbb{R}^d)$. The
operator $\tau_x$ satisfies the following properties (see [18]):
\begin{prop}$ $
\begin{itemize}
\item[i)] For $f\in A_k(\mathbb{R}^d)$ and g a bounded function in $
L^1_k(\mathbb{R}^d)$, we have
\begin{eqnarray*}\int_{\mathbb{R}^d}\tau_x(f)(y)g(y)d\nu_k(y)=\int_{\mathbb{R}^d}f(y)\tau_{-x}(g)(y)d\nu_k(y).\end{eqnarray*}
\item[ii)] $\tau_x(f)(y)=\tau_{-y}(f)(-x)$, \;$x,y\in\mathbb{R}^d$.
\item[iii)] For a radial function $f$ in $L^1_k(\mathbb{R}^d)$, we
have
\begin{eqnarray*}\int_{\mathbb{R}^d}\tau_x(f)(y) d\nu_k(y)=\int_{\mathbb{R}^d}f(y)d\nu_k(y).\end{eqnarray*}
\end{itemize}
\end{prop}
According to ([18], Theorem 3.7), the operator $\tau_x$ can be
extended to the space of radial functions
$L^p_k(\mathbb{R}^d)^{rad},$ $1 \leq p \leq 2$ and we have for a
function $f$ in $L^p_k(\mathbb{R}^d)^{^{rad}}$,
\begin{eqnarray*}\|\tau_x(f)\|_{p,k} \leq \|f\|_{p,k}.\end{eqnarray*}We remark that it is still an open problem whether $\tau_x(f)$ can be defined for all
$f$ in $L^1_k(\mathbb{R}^d)$.

 It was shown in [16], that if $f$ is a radial
function in $S(\mathbb{R}^{d})$ with $f(y)=\widetilde{f}(\|y\|)$,
then
\begin{eqnarray}
  \tau_x (f)(y)=\int_{\mathbb{R}^{d}}\widetilde{f}(A(x,y,\eta))d\mu_{x}(\eta)
\end{eqnarray}where $A(x,y,\eta)=\sqrt{\|x\|^{2}+\|y\|^{2}-2<y,\eta>}$ and $\mu_{x}$ is a probability measure supported in the convex hull
$co(G.x)$ of the $G$-orbit of $x$ in $\mathbb{R}^d$. We observe
that,
\begin{eqnarray}
 \eta \in co(G.x)\Longrightarrow \min_{g\in G}\|g.x+y\|\leq  A(x,y,\eta) \leq \max_{g\in
  G}\|g.x+y\|.
\end{eqnarray}
\section{Behavior at infinity for the Riesz potentials associated to the Dunkl operators}
We study in this section, the behavior at infinity of the Riesz
potential of a non compactly supported function. For
$0<\alpha<2\gamma+d$ and $f\in S(\mathbb{R}^{d})$, we recall that
 the Riesz potential $I_{\alpha}^{k}f$ of a function $f$ is given by
\begin{eqnarray*}
I_{\alpha}^{k}f(x)=2^{\gamma+\frac{d}{2}-\alpha}\frac{\Gamma(\gamma+\frac{d-\alpha}{2})}{\Gamma(\frac{\alpha}{2})}
\int_{\mathbb{R}^{d}} \frac{\tau_{y}f(x)}{\|y\|^{2\gamma+d-\alpha}}
d\nu_{k}(y).
\end{eqnarray*} Throughout this section, for $m>0$ and
$x\in\mathbb{R}^{d}$, we denote  by $V_{x}^{m}$ and $W_{x}^{m}$, the
following sets:
\begin{eqnarray*}
  V_{x}^{m} = \{y\in \mathbb{R}^{d};\,\min_{g\in G}\|g.x+y\|<m\|x\|\},\;\;\;
  W_{x}^{m} = \mathbb{R}^{d}\setminus V_{x}^{m}.
\end{eqnarray*}
\begin{prop} Let $0<\alpha<2\gamma+d$ and $f\in A_k(\mathbb{R}^d)$. Then the function   \begin{eqnarray}
 \Phi: x \longmapsto \int_{0}^{+\infty}s^{\gamma+\frac{d-\alpha}{2}-1}\int_{\mathbb{R}^{d}}f(y)
 \tau_{x}(e^{-s\|.\|^2})(y) d\nu_{k}(y)ds \end{eqnarray}is bounded
 on $\mathbb{R}^d$ and we have
 \begin{eqnarray*}I_{\alpha}^{k}f(x)=\displaystyle\frac{2^{\gamma+\frac{d}{2}-\alpha}}{\Gamma(\frac{\alpha}{2})}\Phi(x),\quad x\in\mathbb{R}^{d}.\end{eqnarray*}
\end{prop}
\begin{pf} For each $x \in \mathbb{R}^{d}$ and $s>0$, we have from
(2.3)\begin{eqnarray*}\tau_{x}(e^{-s\|.\|^2})(y)=\int_{\mathbb{R}^{d}}e^{-s(A(x,y,\eta))^{2}}d\mu_{x}(\eta).\end{eqnarray*}
Since $e^{-s(A(x,y,\eta))^{2}} \leq 1$ and $\mu_{x}$ is a
probability measure, we deduce that
\begin{eqnarray}\tau_{x}(e^{-s\|.\|^2})(y)\leq 1.\end{eqnarray} On the other
hand, using Proposition 2.1, we obtain
\begin{eqnarray}\int_{\mathbb{R}^{d}}\tau_{x}(e^{-s\|.\|^2})(y)  d\nu_{k}(y)= \int_{\mathbb{R}^{d}} e^{-s\|y\|^2} d\nu_{k}(y)= c_k^{-1}(2s)^{-2\gamma-\frac{d}{2}}, \end{eqnarray}
where $c_k$ is the Mehta-type constant (see section 2). Now, for
$f\in A_k(\mathbb{R}^d)$, let us decompose $\Phi(x)$ as a sum of two
terms: $\Phi(x)=\Phi_1(x)+\Phi_2(x)$ where
 \begin{eqnarray*}
 \Phi_1(x)&=&\int_{0}^{1}s^{\gamma+\frac{d-\alpha}{2}-1}\int_{\mathbb{R}^{d}}f(y)\tau_{x}(e^{-s\|.\|^2})(y)
  d\nu_{k}(y)ds, \\
 \Phi_2(x)&=&\int_{1}^{+\infty}s^{\gamma+\frac{d-\alpha}{2}-1}\int_{\mathbb{R}^{d}}f(y)\tau_{x}(e^{-s\|.\|^2})(y)
  d\nu_{k}(y)ds.\end{eqnarray*}
 Using (3.2), we obtain
 \begin{eqnarray}|\Phi_1(x)|\leq\int_{0}^{1}s^{\gamma+\frac{d-\alpha}{2}-1}ds\int_{\mathbb{R}^{d}}|f(y)|
  d\nu_{k}(y)< +\infty.\end{eqnarray}By the fact that $f \in L^{\infty}$ and from (3.3), we
get\begin{eqnarray}|\Phi_2(x)|\leq
c\,\|f\|_{\infty,k}\int_{1}^{+\infty}s^{-\gamma-\frac{\alpha}{2}-1}ds<
+\infty.\end{eqnarray}
 Combining (3.4) and (3.5), we conclude that $\Phi$ is bounded.

 To prove
 $I_{\alpha}^{k}f(x)=\displaystyle\frac{2^{\gamma+\frac{d}{2}-\alpha}}{\Gamma(\frac{\alpha}{2})}\Phi(x),\;
 x\in\mathbb{R}^{d},$ we can see by Proposition 2.1 that
 \begin{eqnarray*}
\int_{\mathbb{R}^{d}}f(y)
 \tau_{x}(e^{-s\|.\|^2})(y) d\nu_{k}(y)=\int_{\mathbb{R}^{d}}\tau_{y}f(x)
  e^{-s\|y\|^2} d\nu_{k}(y).
 \end{eqnarray*} Applying the formula (see [19]),\begin{eqnarray*} \|y\|^{-(2\gamma+d-\alpha)}=\frac{1}{\Gamma(\gamma+\frac{d-\alpha}{2})}
\int_{0}^{+\infty}s^{\gamma+\frac{d-\alpha}{2}-1}
  e^{-s\|y\|^2} ds,\end{eqnarray*} and changing the order of
  integrals in (3.1), it yields
   \begin{eqnarray*}I_{\alpha}^{k}f(x)=\displaystyle\frac{2^{\gamma+\frac{d}{2}-\alpha}}{\Gamma(\frac{\alpha}{2})}\Phi(x),\quad x\in\mathbb{R}^{d}.\end{eqnarray*}
According to (2.3), one can see that
\begin{eqnarray*}
  I_{\alpha}^{k}f(x)=c_{\alpha,k} \int_0^{+\infty}s^{\gamma+\frac{d-\alpha}{2}}\int_{\mathbb{R}^{d}}
f(y)\Big(\int_{\mathbb{R}^{d}}e^{-s(A(x,y,\eta))^{2}}d\mu_{x}(\eta)\Big)d\nu_{k}(y)\frac{ds}{s}.
\end{eqnarray*}
  The proof is complete.
\end{pf}
\begin{thm} Let $0<\alpha<2\gamma+d$ and $f\in A_k(\mathbb{R}^d)$ with $f\geq 0$. For $x\in\mathbb{R}^{d}$, put
 \begin{eqnarray*}
\Psi(x)=\int_{0}^{+\infty}s^{\gamma+\frac{d-\alpha}{2}-1}\int_{\mathbb{R}^{d}}f(y)\tau_{x}(e^{-s\|.\|^2})(y)(1+\|y\|)^{2\gamma+d-\alpha}
 d\nu_{k}(y)ds. \end{eqnarray*} Then $I_{\alpha}^{k}f(x)= O(\|x\|^{\alpha-(2\gamma+d)})$ as $\|x\|\rightarrow +\infty$ if and only if $\Psi$ is bounded.
\end{thm}
\begin{pf}
1) We begin with the sufficiency part. Assume that $\Psi$ is bounded
on $\mathbb{R}^{d}$. Using Proposition 3.1, we can write
 \begin{eqnarray}
  I_{\alpha}^{k}f(x) = c_{\alpha,k}\,\int_{0}^{+\infty}s^{\gamma+\frac{d-\alpha}{2}-1}(I_{1}(x,s)+I_{2}(x,s))ds.
   \end{eqnarray}
   where \begin{eqnarray*}
I_{1}(x,s)&=& \int_{V_{x}^{m}}
f(y)\tau_{x}(e^{-s\|.\|^2})(y)d\nu_{k}(y),\\
I_{2}(x,s) &=& \int_{W_{x}^{m}}
f(y)\tau_{x}(e^{-s\|.\|^2})(y)d\nu_{k}(y).
         \end{eqnarray*}
First, we observe that for $y\in V_{x}^{m}$ and $0<m<1$,
\begin{eqnarray*}
                                               |\|x\|-\|y\|| \leq \min_{g\in G}\|g.x+y\|<m\|x\|,
                                             \end{eqnarray*}
                                             then
\begin{eqnarray}
(1-m)\|x\|<\|y\|.
\end{eqnarray}Since $1\leq \Big(\frac{\|y\|}{1+\|y\|}\Big)^{\alpha-(2\gamma+d)}$, we can assert that
 \begin{eqnarray*}
I_{1}(x,s)\leq \int_{V_{x}^{m}}
f(y)\Big(\frac{\|y\|}{1+\|y\|}\Big)^{\alpha-(2\gamma+d)}\tau_{x}(e^{-s\|.\|^2})(y)d\nu_{k}(y)
\end{eqnarray*}then we obtain by (3.7)
\begin{eqnarray*}
 I_{1}(x,s)\leq c\,\|x\|^{\alpha-(2\gamma+d)}\int_{V_{x}^{m}}f(y)\tau_{x}(e^{-s\|.\|^2})(y)(1+\|y\|)^{(2\gamma+d-\alpha)}d\nu_{k}(y).\end{eqnarray*}
  This yields \begin{eqnarray}\int_{0}^{+\infty}s^{\gamma+\frac{d-\alpha}{2}-1} I_{1}(x,s)ds &\leq & c\,\|x\|^{\alpha-(2\gamma+d)}\Psi(x)\nonumber\\
  &\leq& c\,\|x\|^{\alpha-(2\gamma+d)}.
   \end{eqnarray}
Second, for $s>0$, $x\in\mathbb{R}^{d}$ and $y\in W_{x}^{m}$, we get
by (2.4)
\begin{eqnarray*}
 -s (A(x,y,\eta))^{2}\leq-s(\min_{g\in
 G}\|g.x+y\|)^{2}\leq-s(m\|x\|)^{2},
 \end{eqnarray*} then we obtain from (2.3),
\begin{eqnarray*}
   I_2(x,s)&\leq&\int_{W_{x}^{m}}
f(y)\Big(\int_{\mathbb{R}^{d}}e^{-s(m\|x\|)^{2}}d\mu_{x}(\eta)\Big)d\nu_{k}(y)  \\
  &\leq& \int_{W_{x}^{m}}
f(y)e^{-s(m\|x\|)^{2}}d\nu_{k}(y)\\
&\leq&e^{-s(m\|x\|)^{2}}\|f\|_{1,k},
  \end{eqnarray*}which gives that
\begin{eqnarray*} \int_{0}^{+\infty}s^{\gamma+\frac{d-\alpha}{2}-1}I_{2}(x,s)ds\leq
\|f\|_{1,k}\int_{0}^{+\infty}s^{\gamma+\frac{d-\alpha}{2}-1}e^{-s(m\|x\|)^{2}}ds.
\end{eqnarray*}By the change of variables $t=s(m\|x\|)^{2}$ we get
\begin{eqnarray} \int_{0}^{+\infty}s^{\gamma+\frac{d-\alpha}{2}-1}I_{2}(x,s)ds&\leq& \Gamma(\gamma+\frac{d-\alpha}{2})\;\|f\|_{1,k}(m\|x\|)^{\alpha-(2\gamma+d)}\nonumber\\
&\leq& c\;\|x\|^{\alpha-(2\gamma+d)}
\end{eqnarray}
 Hence from (3.6),
(3.8) and (3.9), we deduce that
$$I_{\alpha}^{k}f(x)=O(\|x\|^{\alpha-(2\gamma+d)})\quad \mbox{as}\quad
\|x\|\rightarrow +\infty.$$ 2) For the necessity part, suppose
$I_{\alpha}^{k}f(x)= O(\|x\|^{\alpha-(2\gamma+d)})$ as
$\|x\|\rightarrow +\infty $. Since for $y\in \mathbb{R}^{d},$
$$(1+\|y\|)^{2\gamma+d-\alpha}\leq
2^{2\gamma+d-\alpha}(1+\|y\|^{2\gamma+d-\alpha}),$$ we can write
from Proposition 3.1 that
\begin{eqnarray}
          \Psi(x)\leq c\,(I_{\alpha}^{k}f(x)+J_{1}+J_{2}),
      \end{eqnarray}where
      \begin{eqnarray*}
 J_{1}&=&
 \int_{0}^{+\infty}s^{\gamma+\frac{d-\alpha}{2}-1}\int_{W_{x}^{m}}f(y)\tau_{x}(e^{-s\|.\|^2})(y)\|y\|^{2\gamma+d-\alpha}d\nu_{k}(y)ds,\\
 J_{2}&=& \int_{0}^{+\infty}s^{\gamma+\frac{d-\alpha}{2}-1}\int_{V_{x}^{m}}f(y)\tau_{x}(e^{-s\|.\|^2})(y)\|y\|^{2\gamma+d-\alpha}d\nu_{k}(y)ds.\end{eqnarray*}
From (2.4) and for $y\in W_{x}^{m}$, we have $\displaystyle
A(x,y,\eta)\geq \min_{g\in G}\|g.x+y\|\geq m\|x\|$. Then by the fact
that $\mu_{x}$ is a probability measure and using (2.3) and Fubini's
theorem, we obtain
\begin{eqnarray*}
   J_{1}\leq \int_{W_{x}^{m}}f(y)\|y\|^{2\gamma+d-\alpha}\int_{0}^{+\infty}s^{\gamma+\frac{d-\alpha}{2}-1} \displaystyle{e^{-s(\min_{g\in G}\|g.x+y\|)^{2}}}  d\nu_{k}(y)ds.\end{eqnarray*}
By the change of variables $\displaystyle t=s(\min_{g\in
G}\|g.x+y\|)^{2}$, we can assert that
\begin{eqnarray*}
 J_{1}\leq \Gamma(\gamma+\frac{d-\alpha}{2})\int_{W_{x}^{m}}f(y)
(\displaystyle\min_{g\in
G}\|g.x+y\|)^{\alpha-2\gamma-d}\|y\|^{2\gamma+d-\alpha}
d\nu_{k}(y)ds.
\end{eqnarray*}
Since $\|y\|\leq \|x\|+\displaystyle\min_{g\in G}\|g.x+y\|$, this
gives
\begin{eqnarray*}
   \min_{g\in G}\|g.x+y\|^{\alpha-(2\gamma+d)}\|y\|^{2\gamma+d-\alpha}
  &=& \Big(\frac{\|y\|}{\displaystyle\min_{g\in G}\|g.x+y\|}\Big)^{2\gamma+d-\alpha}\\
&\leq & \Big(1+\frac{\|x\|}{\displaystyle\min_{g\in G}\|g.x+y\|}\Big)^{2\gamma+d-\alpha}\\
  &\leq& \Big(1+\frac{1}{m}\Big)^{2\gamma+d-\alpha},
  \end{eqnarray*}hence, we deduce that
  \begin{eqnarray}
   J_{1}&\leq&\Gamma(\gamma+\frac{d-\alpha}{2}) (1+\frac{1}{m})^{2\gamma+d-\alpha}\int_{W_{x}^{m}}f(y)d\nu_{k}(y)\nonumber \\
     &\leq& c\,\|f\|_{1,k}.
  \end{eqnarray}
  If $y\in V_{x}^{m}$, then $\|y\|\leq (1+m)\|x\|$. This yields from Proposition 3.1 that\begin{eqnarray}
       J_{2}&\leq& c\,\|x\|^{2\gamma+d-\alpha} \int_{0}^{+\infty}s^{\gamma+\frac{d-\alpha}{2}-1}
       \int_{V_{x}^{m}}f(y)\tau_{x}(e^{-s\|.\|^2})(y)d\nu_{k}(y)ds.\nonumber
       \\
        &\leq& c\,\|x\|^{2\gamma+d-\alpha}I_{\alpha}^{k}f(x)<+\infty.
      \end{eqnarray}
      Using (3.10), (3.11) and (3.12), we conclude that $\Psi$ is
      bounded. This completes the proof of the theorem.
            \end{pf}
 \begin{rem} Take $f$ in $S(\mathbb{R}^{d})$, then $f\in A_k(\mathbb{R}^d)$. We
shall prove that the function $\Psi$ given by\begin{eqnarray*}
 \Psi(x)= \int_{0}^{+\infty}s^{\gamma+\frac{d-\alpha}{2}-1}\int_{\mathbb{R}^{d}}|f(y)|\,\tau_{x}(e^{-s\|.\|^2})(y)
 (1+\|y\|)^{2\gamma+d-\alpha}d\nu_{k}(y)ds\end{eqnarray*} is bounded on $\mathbb{R}^{d}$.
 We write $\Psi(x)=\Psi_1(x)+\Psi_2(x)$ where
  \begin{eqnarray*}
  \Psi_1(x)&=&\int_{0}^{1}s^{\gamma+\frac{d-\alpha}{2}-1}\int_{\mathbb{R}^{d}}|f(y)|\,\tau_{x}(e^{-s\|.\|^2})(y)
 (1+\|y\|)^{2\gamma+d-\alpha}d\nu_{k}(y)ds,\\
 \Psi_2(x)&=&\int_{1}^{+\infty}s^{\gamma+\frac{d-\alpha}{2}-1}\int_{\mathbb{R}^{d}}|f(y)|\,\tau_{x}(e^{-s\|.\|^2})(y)
 (1+\|y\|)^{2\gamma+d-\alpha}d\nu_{k}(y)ds.\end{eqnarray*}
 Since $f\in S(\mathbb{R}^{d})$, we obtain using (3.2)
 \begin{eqnarray}\Psi_1(x)\leq\int_{0}^{1}s^{\gamma+\frac{d-\alpha}{2}-1}ds\int_{\mathbb{R}^{d}}|f(y)|
 (1+\|y\|)^{2\gamma+d-\alpha}d\nu_{k}(y)< +\infty.\end{eqnarray}By the fact that $y\rightarrow
 (1+\|y\|)^{2\gamma+d-\alpha}f(y)$ is bounded on $\mathbb{R}^{d}$ and from (3.3), we get
 \begin{eqnarray}\Psi_2(x)\leq c\int_{1}^{+\infty}s^{-\gamma-\frac{\alpha}{2}-1}ds < +\infty.\end{eqnarray}
 Combining (3.13) and (3.14), we conclude that $\Psi$ is bounded. \\Since $|I_{\alpha}^{k}f(x)|\leq I_{\alpha}^{k}|f|(x),$
 $x \in \mathbb{R}^{d}$, we deduce from Theorem 3.1 that
 $$\; I_{\alpha}^{k}f(x)=O(\|x\|^{\alpha-(2\gamma+d)})\quad \mbox{as}\quad  \|x\| \rightarrow +\infty.$$
\end{rem}
\begin{ex} Take the function, $f(x)=e^{-\|x\|}$, $x\in\mathbb{R}^{d}$. It was shown in [18] that
 $\displaystyle\mathcal{F}_{k}(f)(x)=c_{d,k}\frac{1}{(1+\|x\|^{2})^{\gamma+\frac{d+1}{2}}},$ where
      $ c_{d,k}=\displaystyle 2^{\gamma+\frac{d}{2}}\frac{\Gamma(\gamma+\frac{d+1}{2})}{\sqrt{\pi}}$.\\
      We can see that $f\in A_{k}(\mathbb{R}^{d})$ and the function
       $g:x\longmapsto e^{-\|x\|}(1+\|x\|)^{2\gamma+d-\alpha}$ is in $L_{k}^{1}(\mathbb{R}^{d})$ and  bounded.\\
      By proceeding in the same manner as in Remark 3.1, we write,
       \begin{eqnarray*}
 \Psi(x)&=&
 \int_{0}^{+\infty}s^{\gamma+\frac{d-\alpha}{2}-1}\int_{\mathbb{R}^{d}}g(y)\tau_{x}(e^{-s\|.\|^2})(y)d\nu_{k}(y)ds.\\
 &=& \int_{0}^{1}s^{\gamma+\frac{d-\alpha}{2}-1}\int_{\mathbb{R}^{d}}g(y)\tau_{x}(e^{-s\|.\|^2})(y)d\nu_{k}(y)ds. \\
 &+& \int_{1}^{+\infty}s^{\gamma+\frac{d-\alpha}{2}-1}\int_{\mathbb{R}^{d}}g(y)\tau_{x}(e^{-s\|.\|^2})(y)d\nu_{k}(y)ds. \\
 &=&\Psi_{1}(x)+\Psi_{2}(x).
                        \end{eqnarray*}
 Using respectively (3.2) for $\Psi_1$ and (3.3) for $\Psi_2$, we can
 assert $\Psi_1$ and $\Psi_2$ are bounded. We conclude that $\Psi$ is bounded and applying Theorem 3.1, we
 obtain that $I_{\alpha}^{k}f(x)= O(\|x\|^{\alpha-(2\gamma+d)})$ as
$\|x\|\rightarrow +\infty$.
\end{ex}
\section{Weighted norm inequalities}
In this section, we prove for the Riesz potential $I_{\alpha}^{k}$,
weighted norm inequalities with sufficient conditions on
non-negative pairs of weight functions. We denote by $p'$ the
conjugate of $p$ for $1<p<+\infty$. The proof requires a useful
well-known facts and results which we shall now state in the
following.
\begin{rem}$ $\\
1/ (see [6]) (Hardy inequalities) If $\mu$ and $\vartheta$ are
locally integrable
 weight functions on $(0,+\infty)$ and $1<p\leq q<+\infty$, then there is a constant $c>0$ such that for
all non-negative Lebesgue measurable function $f$ on $(0,+\infty)$,
the inequality
\begin{eqnarray}\Big(\int_{0}^{+\infty}\Big[\int_{0}^{t}f(s)ds\Big]^{q}\mu(t)dt\Big)^{\frac{1}{q}}\leq
 c\, \Big(\int_{0}^{+\infty}(f(t))^{p}\vartheta(t)dt\Big)^{\frac{1}{p}}\end{eqnarray}
 is satisfied if and only if
\begin{eqnarray}\displaystyle\sup_{s>0}\Big(\int_{s}^{+\infty}\mu(t)dt\Big)^{\frac{1}{q}}
\Big(\int_{0}^{s}(\vartheta(t))^{1-p'}dt\Big)^{\frac{1}{p'}}<+
\infty.
\end{eqnarray}
Similarly for the dual operator,
\begin{eqnarray}\Big(\int_{0}^{+\infty}\Big[\int_{t}^{+\infty}f(s)ds\Big]^{q}\mu(t)dt\Big)^{\frac{1}{q}}\leq
 c\Big(\int_{0}^{+\infty}(f(t))^{p}\vartheta(t)dt\Big)^{\frac{1}{p}}
 \end{eqnarray}
 is satisfied if and only if
\begin{eqnarray}\displaystyle\sup_{s>0}\Big(\int_{0}^{s}\mu(t)dt\Big)^{\frac{1}{q}}
\Big(\int_{s}^{+\infty}(\vartheta(t))^{1-p'}dt\Big)^{\frac{1}{p'}}<+
\infty.
\end{eqnarray}
2/ Let $f$ be a complex-valued $\nu_k$-measurable function on
$\mathbb{R}^{d}$. The distribution function $D_f$ of $f$ is defined
for all $s\geq0$ by
$$D_{f}(s)=\nu_k(\{x\in\mathbb{R}^{d}\,:\; |f(x)|>s\}).$$ The decreasing
rearrangement of $f$ is the function $f^*$ given for all $t\geq0$ by
$$f^{*}(t)=inf\{s\geq0 \,:\; D_{f}(s)\leq t\}.$$
We list some known results:
\begin{itemize}
\item[$\bullet$] Let $f\in
L^p_{k}(\mathbb{R}^d)$ and $1\leq p<+\infty$, then
\begin{eqnarray*}\int_{\mathbb{R}^{d}}|f(x)|^{p}
d\nu_{k}(x)=p\int_{0}^{+\infty}s^{p-1}D_{f}(s)ds=\int_{0}^{+\infty}(f^{*}(t))^{p}dt.\end{eqnarray*}
\item[$\bullet$] (see [11]) (Hardy-Littlewood rearrangement inequality)\\
Let $f$ and $\upsilon$ be non negative $\nu_{k}$-measurable
functions on $\mathbb{R}^{d}$, then
\begin{eqnarray}\int_{\mathbb{R}^{d}}f(x)\upsilon(x)d\nu_{k}(x)\leq\int_{0}^{+\infty}f^{*}(t)\upsilon^{*}(t)dt
\end{eqnarray}
and
\begin{eqnarray}\int_{0}^{+\infty}f^{*}(t)\frac{1}
{(\frac{1}{\upsilon})^{*}(t)}dt\leq\int_{\mathbb{R}^{d}}f(x)\upsilon(x)d\nu_{k}(x).
\end{eqnarray}
\item[$\bullet$] (see A. P. Calder\'{o}n [7]) Let $1\leq p_{1}< p_{2}<\infty
$ and $1 \leq q_{1}<q_{2}<+ \infty $. A linear operator
$\mathcal{L}$ satisfies the weak-type hypotheses $(p_{1},q_{1})$ and
$(p_{2},q_{2})$  if and only if
\begin{eqnarray}
  (\mathcal{L}f)^{\ast}(t)&\leq& c\,\Big(t^{-\frac{1}{q_{1}}}\int_{0}^{ t^{\frac{\lambda_{1}}{\lambda_{2}}}}s^{\frac{1}{p_{1}}-1}f^{\ast}(s)ds+
   t^{-\frac{1}{q_{2}}}\int_{t^{\frac{\lambda_{1}}{\lambda_{2}}}}^{+\infty}s^{\frac{1}{p_{2}}-1}f^{\ast}(s)ds
   \Big),\nonumber\\ &&
\end{eqnarray} where $\lambda_{1} =\frac{1}{q_{1}}-\frac{1}{q_{2}}$
and $\lambda_{2}=\frac{1}{p_{1}}-\frac{1}{p_{2}}$.
\end{itemize}
\end{rem}
\begin{ex} Let $\delta<0$ and $\beta>0$. Take $u(x)=\|x\|^{\delta}$, $v(x)=\|x\|^{\beta}$,
 and $f(x)=\chi_{(0,r)}(\|x\|)$, $x\in\mathbb{R}^{d}$, $r>0$. Then using
(2.1), we
have for $s\geq0$, \begin{eqnarray*}D_{u}(s)&=&\nu_{k}\Big(\{x\in\mathbb{R}^{d}\,:\;\|x\|^{\delta}>s\}\Big)\\
&=&\nu_{k}\Big(B(0,s^{\frac{1}{\delta}})\Big)=
\frac{d_{k}}{2\gamma+d}\;s^{\frac{2\gamma+d}{\delta}},\end{eqnarray*}
\begin{eqnarray*}D_{\frac{1}{v}}(s)&=&\nu_{k}\Big(\{x\in\mathbb{R}^{d}\,:\;\|x\|^{-\beta}>s\}\Big)\\
&=&\nu_{k}\Big(B(0,s^{-\frac{1}{\beta}})\Big)=\frac{d_{k}}{2\gamma+d}\;s^{-\frac{2\gamma+d}{\beta}},\end{eqnarray*}
and \begin{eqnarray*}D_{f}(s)&=&\nu_{k}\Big(\{x\in\mathbb{R}^{d}\,:\;\chi_{(0,r)}(\|x\|)>s\}\Big)\\
&=& \nu_{k}\Big(B(0,1)\Big)\,r^{2\gamma+d}\,\chi_{(0,1)}(s)\\
&=&\frac{d_k}{2\gamma+d}\,r^{2\gamma+d}\,\chi_{(0,1)}(s).\end{eqnarray*}
Note that, $\displaystyle
d_k=\frac{c^{-1}_k}{2^{\gamma+\frac{d}{2} -1} \Gamma(\gamma +
\frac{d}{2})},$ this yields
$$\frac{d_k}{2\gamma+d}=  \frac{c^{-1}_k}{2^{\gamma +\frac{d}{2}}
\Gamma(\gamma + \frac{d}{2}+1)}.$$ This gives for $t\geq0$,
\begin{eqnarray*}u^{*}(t)=inf\{s\geq 0\,:\; D_{u}(s)\leq t\}=
\Big(\frac{2\gamma+d}{d_{k}}\Big)^{\frac{\delta}{2\gamma+d}}\;t^{\frac{\delta}{2\gamma+d}},\end{eqnarray*}
\begin{eqnarray*}(\frac{1}{v})^{*}(t)=inf\{s\geq 0\,:\; D_{\frac{1}{v}}(s)\leq t\}
=\Big(\frac{2\gamma+d}{d_{k}}\Big)^{-\frac{\beta}{2\gamma+d}}\;t^{-\frac{\beta}{2\gamma+d}},\end{eqnarray*}
and \begin{eqnarray*}f^*(t)=\chi_{(0,R)}(t)\;\;\mbox{where}\;\;
R=\frac{d_k}{2\gamma+d}\,r^{2\gamma+d}.\end{eqnarray*} Hence, using
(2.1) again, we obtain for $-(2\gamma+d)<\delta$
\begin{eqnarray*}\int_{\mathbb{R}^{d}}
f(x)u(x)d\nu_{k}(x)&=& \frac{d_k
}{\delta+2\gamma+d}\;r^{\delta+2\gamma+d}\\&=&\int_{0}^{+\infty}f^{*}(t)u^{*}(t)dt,\end{eqnarray*}
and
\begin{eqnarray*}\int_{0}^{+\infty}f^{*}(t)\frac{1}
{(\frac{1}{v})^{*}(t)}dt&=& \frac{d_k
}{\beta+2\gamma+d}\;r^{\beta+2\gamma+d}\\&=&\int_{\mathbb{R}^{d}}
f(x)v(x)d\nu_{k}(x),\end{eqnarray*} giving equalities for (4.5) and
(4.6) in these cases.
\end{ex}
In the following theorem, we prove weighted norm inequalities for
the Riesz potential $I_{\alpha}^{k}$.
\begin{thm} Let $0<\alpha<2\gamma+d$, $1<r<\frac{2\gamma+d}{\alpha}$ and $u$,$v$ be non negative $\nu_{k}$-locally integrable weight functions on
$\mathbb{R}^{d}$. Then for $1< p \leq q <+ \infty $,
$I_{\alpha}^{k}$ can be extended to a bounded operator from
$L_{k,v}^{p}(\mathbb{R}^{d})$
  to $L_{k,u}^{q}(\mathbb{R}^{d})$ and the inequality \begin{eqnarray*}\|I_{\alpha}^{k}f\|_{q,k,u}\leq c\,\|f\|_{p,k,v}\end{eqnarray*}
   holds with the following conditions on $u$ and
  $v$:
\begin{eqnarray}
  \sup_{s>0}\Big(\int_{s}^{+\infty}u^{\ast}(t)t^{-q(1-\frac{\alpha}{2\gamma+d})}dt
  \Big)^{\frac{1}{q}} \Big(\int_{0}^{s}[(\frac{1}{v})^{\ast}(t)]^{(p'-1)}dt\Big)^{\frac{1}{p'}}<+\infty
  \end{eqnarray}
   and
   \begin{eqnarray}
   \sup_{s>0} \Big(\int_{0}^{s}u^{\ast}(t)t^{-q(\frac{1}{r}-\frac{\alpha}{2\gamma+d})}dt\Big)^{\frac{1}{q}}
\Big(\int_{s}^{+\infty}[(\frac{1}{v})^{\ast}(t)^{(p'-1)}t^{p'(\frac{1}{r}-1)}dt\Big)^{\frac{1}{p'}}
<+\infty.
\end{eqnarray}
\end{thm}
\begin{pf}
It was shown in [12] that: $f\rightarrow I_{\alpha}^{k}$ can be
extended to a mapping of weak-type
$(p_1,q_1)=(1,\frac{1}{1-\frac{\alpha}{2\gamma+d}})$ and a bounded
operator from $L^{p_2}_{k}(\mathbb{R}^d)$ to
$L^{q_2}_{k}(\mathbb{R}^d)$ with
$(p_2,q_2)=(r,\frac{1}{\frac{1}{r}-\frac{\alpha}{2\gamma+d}})$. Put
$\lambda_{1} =\lambda_{2}=1-\frac{1}{r}$. Let
$f\in\mathcal{S}(\mathbb{R}^{d})$. Using (4.7) and applying
Minkowski's inequality, we have\\
$\displaystyle\Big(\int_{0}^{+\infty}[(I_{\alpha}^{k}f)^{\ast}(t)]^{q}u^{\ast}(t)dt\Big)^{\frac{1}{q}}$
\begin{eqnarray*}
 &\leq & c\,
   \Big[\int_{0}^{+\infty}u^{\ast}(t)t^{-\frac{q}{q_{1}}}
  \Big(\int_{0}^{t^{\frac{\lambda_{1}}{\lambda_{2}}}}s^{\frac{1}{p_{1}}-1}f^{\ast}(s)ds\Big)^{q}dt\Big]^{\frac{1}{q}} \\
   && +\,c\,\Big[\int_{0}^{+\infty}u^{\ast}(t)t^{-\frac{q}{q_{2}}}
  \Big(\int_{t^{\frac{\lambda_{1}}{\lambda_{2}}}}^{+\infty}s^{\frac{1}{p_{2}}-1}f^{\ast}(s)ds\Big)^{q}dt\Big]^{\frac{1}{q}}.
\end{eqnarray*}
By means of change of variable in the right side, we obtain\\
$\displaystyle\Big(\int_{0}^{+\infty}[(I_{\alpha}^{k}f)^{\ast}(t)]^{q}u^{\ast}(t)dt\Big)^{\frac{1}{q}}$
\begin{eqnarray}
   &\leq & c\,
 \Big[\int_{0}^{+\infty}u^{\ast}(t^{\frac{\lambda_{2}}{\lambda_{1}}})t^{\frac{\lambda_{2}}{\lambda_{1}}(1-\frac{q}{q_{1}})-1}
  [\int_{0}^{t}s^{\frac{1}{p_{1}}-1}f^{\ast}(s)ds]^{q}dt\Big]^{\frac{1}{q}}\nonumber\\
  &&+\,
  c\,\Big[\int_{0}^{+\infty}u^{\ast}(t^{\frac{\lambda_{2}}{\lambda_{1}}})t^{\frac{\lambda_{2}}{\lambda_{1}}(1-\frac{q}{q_{2}})-1}
  [\int_{t}^{+\infty}s^{\frac{1}{p_{2}}-1}f^{\ast}(s)ds]^{q}dt\Big]^{\frac{1}{q}}\nonumber \\
 & =& I_{1}+I_{2}.
 \end{eqnarray}
Applying (4.1) and (4.2) for $I_1$, we can assert that
 \begin{eqnarray} I_{1}\leq \Big(\int_{0}^{+\infty}[(\frac{1}{v})^{\ast}(t)]^{-1}[f^{\ast}(t)]^{p}dt\Big)^{\frac{1}{p}}\end{eqnarray}
 if and only if
 \begin{eqnarray*}
   \sup_{s>0}\Big(\int_{s}^{+\infty}u^{\ast}(t^{\frac{\lambda_{2}}{\lambda_{1}}})t^{\frac{\lambda_{2}}{\lambda_{1}}(1-\frac{q}{q_{1}})-1}dt\Big)^{\frac{1}{q}}
    \Big(\int_{0}^{s}[(\frac{1}{v})^{\ast}(t)]^{(p'-1)}t^{p'(\frac{1}{p_{1}}-1)}dt\Big)^{\frac{1}{p'}}\leq
    +\infty.
 \end{eqnarray*}
Then if we replace $s$ by $s^{\frac{1}{\lambda_{2}}}$ in this
condition, it's easy to see that if we use a change of variable in
the first integral of the expression, we obtain (4.8).  \\
 Similarly by applying (4.3) and (4.4) for $I_2$, we get
\begin{eqnarray}
    I_{2}\leq \Big(\int_{0}^{+\infty}[(\frac{1}{v})^{\ast}(t)]^{-1}[f^{\ast}(t)]^{p}dt\Big)^{\frac{1}{p}}
 \end{eqnarray}
 if and only if
 \begin{eqnarray*}
    \sup_{s>0}\Big(\int_{0}^{s}u^{\ast}(t^{\frac{\lambda_{2}}{\lambda_{1}}})t^{\frac{\lambda_{2}}{\lambda_{1}}(1-\frac{q}{q_{2}})-1}dt\Big)^{\frac{1}{q}}
    \Big(\int_{s}^{+\infty}[(\frac{1}{v})^{\ast}(t)]^{(p'-1)}t^{p'(\frac{1}{p_{2}}-1)}dt\Big)^{\frac{1}{p'}}\leq
    +\infty,
 \end{eqnarray*}
 which is equivalent to (4.9).\\
 Combining (4.10), (4.11) and (4.12), it yields
 \begin{eqnarray}
 \displaystyle\Big(\int_{0}^{+\infty}[(I_{\alpha}^{k}f)^{\ast}(t)]^{q}u^{\ast}(t)dt\Big)^{\frac{1}{q}}
  \leq
  c\,\Big(\int_{0}^{+\infty}[(\frac{1}{v})^{\ast}(t)]^{-1}[f^{\star}(t)]^{p}dt\Big)^{\frac{1}{p}}.
  \end{eqnarray}
 Using (4.5) on the left side and (4.6) on the right side of (4.13), we
 obtain by density of $\mathcal{S}(\mathbb{R}^{d})$ in
 $L_{k,v}^{p}(\mathbb{R}^{d})$, $1\leq p<+\infty$,
 \begin{eqnarray*}
 \Big(\int_{\mathbb{R}^{d}}[(I_{\alpha}^{k}f)(x)]^{q}u(x)d\nu_{k}(x)\Big)^{\frac{1}{q}}
 \leq c\,
 \Big(\int_{\mathbb{R}^{d}}[f(x)]^{p}v(x)d\nu_{k}(x)\Big)^{\frac{1}{p}}.
  \end{eqnarray*}
This completes the proof.
\end{pf}

As consequence of Theorem 4.1 for power weights, we obtain the
result below.
\begin{cor}
Let $0<\alpha<2\gamma+d$ and $1<p<\frac{2\gamma+d}{\alpha}$. For
$\delta,\beta$ such that $\delta < 0 $, $0< \beta =\delta+\alpha p
<(2\gamma+d)(p-1)$ and $f\in L_{k,v}^{p}(\mathbb{R}^{d})$ with
$v=\|.\|^\beta$, we have
\begin{eqnarray*}
\Big(\int_{\mathbb{R}^{d}}|I_{\alpha}^{k}f(x)|^{p}\|x\|^{\delta}d\nu_{k}(x)\Big)^{\frac{1}{p}}
\leq c\,
  \Big(\int_{\mathbb{R}^{d}}|f(x)|^{p}\|x\|^{\beta}d\nu_{k}(x)\Big)^{\frac{1}{p}}
\end{eqnarray*}
\end{cor}
\begin{pf}
From Example 4.1, we have for $\delta<0$ and
$\beta>0$\begin{eqnarray*}
u^{\ast}(t)=\Big(\frac{2\gamma+d}{d_{k}}\Big)^{\frac{\delta}{2\gamma+d}}t^{\frac{\delta}{2\gamma+d}}
\quad\mbox{and}\quad
(\frac{1}{v})^{\ast}(t)=\Big(\frac{2\gamma+d}{d_{k}}\Big)^{-\frac{\beta}{2\gamma+d}}t^{-\frac{\beta}{2\gamma+d}},\end{eqnarray*}
then if we take $p=q=r$ in Theorem 4.1, the boundedness conditions
(4.8) and (4.9) are valid if and only if
 $$ \left\{\begin{array}{lll}0<\beta<(2\gamma+d)(p-1),\\
&&\\\beta=\delta+\alpha p\,.\end{array}\right.$$ Under these
conditions and from Theorem 4.1, we obtain our result.
\end{pf}
\begin{rem} The boundedness of Riesz potentials can be used to establish the
boundedness properties of the fractional maximal operator given by
\begin{eqnarray*}M_{k,\alpha}f(x)= \displaystyle \sup_{r>0}\frac{1}{m_{k}\,r^{d+2\gamma-\alpha}}
 \int_{\mathbb{R}^{d}}|f(y)|\,\tau_{x}\chi_{B_{r}}(y)d\nu_{k}(y),\quad x\in \mathbb{R}^{d},\end{eqnarray*}
where $$m_{k}=\Big(c_k\,
2^{\gamma+\frac{d}{2}}\Gamma(\gamma+\frac{d}{2}+1)\Big)^{\frac{\alpha}{d+2\gamma}-1}$$
and $\chi_{B_{r}}$ is the characteristic function of the ball
$B_{r}=B(0,r)$. This follows from the fact that
$$M_{k,\alpha}f(x)\leq c\,I_{\alpha}^{k}(|f|)(x),\;x\in \mathbb{R}^{d},$$
which gives for $M_{k,\alpha}$, the same results obtained in Theorem
4.1 and Corollary 4.1.
 \end{rem}

In order to prove a weighted generalized Sobolev inequality, we need
some useful results that we state in the following remark.
\begin{rem}$ $ \\
1/ (see[19]) In Dunkl setting the Riesz transforms  are the
operators $\mathcal{R}_{j}$, $j=1...d,$ defined on
$L^{2}_{k}(\mathbb{R}^{d})$ by
\begin{eqnarray*}
  \mathcal{R}_{j}(f)(x) &=& 2^{\frac{\ell_{k}-1}{2}}\frac{\Gamma(\frac{\ell_{k}}{2})}{\sqrt{\pi}} \,\lim_{\epsilon\rightarrow 0}\int_{\|y\|>\epsilon}\tau_{x}(f)(-y)\frac{y_{j}}{\|y\|^{p_{k}}}d\nu_{k}(y),x\in\mathbb{R}^{d}
\end{eqnarray*}where
\begin{eqnarray*}
 \ell_{k}=2\gamma+d+1.
\end{eqnarray*}
$\bullet$ The Riesz transform $\mathcal{R}_{j}$ is a multiplier
operator with
\begin{eqnarray}
\mathcal{F}_{k}(\mathcal{R}_{j}(f))(\xi) =
\frac{-i\xi_{j}}{\|\xi\|}\,\mathcal{F}_{k}(f)(\xi),\;1\leq j\leq d,
\quad f\in
  \mathcal{S}(\mathbb{R}^{d}).
\end{eqnarray}
$\bullet$ Let $0<\alpha<2\gamma+d $. The identity
\begin{eqnarray}\mathcal{F}_{k}(I_{\alpha}^{k}f)(x)=\|x\|^{-\alpha}\mathcal{F}_{k}(f)(x)\end{eqnarray}
holds in the sense that\begin{eqnarray*}
\int_{\mathbb{R}^{d}}I_{\alpha}^{k}f(x)g(x)d\nu_{k}(x)=\int_{\mathbb{R}^{d}}\mathcal{F}_{k}(f)(x)\|x\|^{-\alpha}\mathcal{F}_{k}(g)(x)d\nu_{k}(x),\end{eqnarray*}
whenever $f,g\in \mathcal{S}(\mathbb{R}^{d})$.\\\\
2/ (see [5]) The Riesz transform $\mathcal{R}_{j}$,$\;1\leq j\leq
d$, can be extended to a bounded operator from
$L_{k}^{p}(\mathbb{R}^{d})$ into it self for $1<p<+\infty$ and we
have
\begin{eqnarray}\|\mathcal{R}_{j}(f)\|_{p,k}\leq c\,\|f\|_{p,k}.\end{eqnarray}
\end{rem}

Now, we give in the following theorem a weighted generalized Sobolev
inequality.
\begin{thm} Let $u$ be a non-negative $\nu_{k}$-locally integrable function on
$\mathbb{R}^{d}$ and $1<r<2\gamma+d$. Then for $1<p\leq q<+\infty$
such that $p<r$ and $f\in \mathcal{D}(\mathbb{R}^{d})$, the
inequality \begin{eqnarray*}\|f\|_{q,k,u} \leq
c\,\|\nabla_{k}f\|_{p,k},\end{eqnarray*} holds with the following
conditions on $u$:
\begin{eqnarray}
\Big(\int_{s}^{\infty}u^{\ast}(t)t^{-q(1-\frac{1}{2\gamma+d})}dt\Big)^{\frac{1}{q}}
\leq c\,s^{\frac{1}{p}-1}
    \end{eqnarray}
   and
   \begin{eqnarray}
    \Big(\int_{0}^{s}u^{\ast}(t)t^{-q(\frac{1}{r}-\frac{1}{2\gamma+d})}dt\Big)^{\frac{1}{q}}
\leq c\,s^{\frac{1}{p}-\frac{1}{r}},
\end{eqnarray} for all $s>0$. Here $\nabla_{k}f=(T_{1}f,...,T_{d}f)$ and
$\displaystyle|\nabla_{k}f|=\Big(\sum_{j=1}^{d}|T_{j}f|^{2}\Big)^{\frac{1}{2}}$.
\end{thm}
\begin{pf}
For $f\in \mathcal{D}(\mathbb{R}^{d})$, we write
\begin{eqnarray*}
  \mathcal{F}_{k}(f)(\xi)=\frac{1}{\|\xi\|}\sum_{j=1}^{d}\frac{-i\xi_{j}}{\|\xi\|}\,i\xi_{j}\,\mathcal{F}_{k}(f)(\xi),
 \end{eqnarray*}then by (2.2) and (4.14), we get
 \begin{eqnarray*}\mathcal{F}_{k}(f)(\xi)&=&\frac{1}{\|\xi\|}\sum_{j=1}^{d}\frac{-i\xi_{j}}{\|\xi\|}\,\mathcal{F}_{k}(T_{j}f)(\xi)\\
 &=&\frac{1}{\|\xi\|}\sum_{j=1}^{d} \mathcal{F}_{k}(\mathcal{R}_{j}(T_{j}f))(\xi)\\
 &=&\frac{1}{\|\xi\|}\mathcal{F}_{k}\Big(\sum_{j=1}^{d} \mathcal{R}_{j}(T_{j}f)\Big)(\xi).
\end{eqnarray*}
This yields from (4.15) that\begin{eqnarray*}
\mathcal{F}_{k}(f)(\xi)=\mathcal{F}_{k}\Big[I_{1}^{k}\Big(\sum_{j=1}^{d}\mathcal{R}_{j}(T_{j}f)\Big)\Big](\xi),
\end{eqnarray*} which gives the following identity,
\begin{eqnarray*}
f=I_{1}^{k}\Big(\sum_{j=1}^{d}\mathcal{R}_{j}(T_{j}f)\Big).
\end{eqnarray*}
Now, observe that the conditions (4.17) and (4.18) are equivalent to
(4.8) and (4.9) with $v\equiv1$ and $\alpha =1$, then using Theorem
4.1, we obtain
\begin{eqnarray*}
  \|f\|_{q,k,u} &=& \|I_{1}^{k}\Big(\sum_{j=1}^{d}\mathcal{R}_{j}(T_{j}f)\Big)\|_{q,k,u}\nonumber \\
   &\leq&c\,\|\mathcal{R}_{j}(\sum_{j=1}^{d}(T_{j}f))\|_{p,k},\end{eqnarray*}
which gives from (4.16) that \begin{eqnarray*}\|f\|_{q,k,u}  &\leq&c\,\|\sum_{j=1}^{d}(T_{j}f)\|_{p,k} \\
   &\leq&c\,\|\nabla_{k}f\|_{p,k}.
\end{eqnarray*}Our result is proved.
\end{pf}
\begin{cor}
Let $1<p<2\gamma+d$ and $1<p\leq q<+\infty$. Then for $\delta < 0$
such that $\delta=q\,[(2\gamma+d)(\frac{1}{p}-\frac{1}{q})-1]$, we
have for $f\in \mathcal{D}(\mathbb{R}^{d})$
\begin{eqnarray*}
\Big(\int_{\mathbb{R}^{d}}|f(x)|^{q}\|x\|^{\delta}d\nu_{k}(x)\Big)^{\frac{1}{q}}
\leq c\,\|\nabla_{k}f\|_{p,k}.
\end{eqnarray*}
\end{cor}
\begin{pf}
For $\delta<0$, if we take $u(x)=\|x\|^\delta$, $x\in \mathbb{R}^{d}
$ in Theorem 4.2, the boundedness conditions (4.17) and (4.18) are
valid if and only if
\begin{eqnarray*}\delta=q\,\Big[(2\gamma+d)\Big(\frac{1}{p}-\frac{1}{q}\Big)-1\Big].\end{eqnarray*} Under this
condition, we obtain our result.
\end{pf}
\begin{rem}
The case $u\equiv1$, $1<p<2\gamma+d$ and
$\frac{1}{q}=\frac{1}{p}-\frac{1}{2\gamma+d}$ was obtained in [5]
and gives the generalized Sobolev inequality
\begin{eqnarray*}
\|f\|_{q,k} \leq c\,\|\nabla_{k}f\|_{p,k}.
\end{eqnarray*}
\end{rem}


\begin{thebibliography}{}
\bibitem{} C. Abdelkefi, J. Ph. Anker, F. Sassi and M. Sifi,
 Besov-type spaces on $\mathbb{R}^d $ and integrability for the Dunkl
transform. Symmetry, Integrability and Geometry: Methods and
Applications, SIGMA 5 (2009), 019, 15 pages.
\bibitem{} C. Abdelkefi, Dunkl
operators on $\mathbb{R}^d$ and uncentered maximal function.
 J. Lie Theory 20 (2010), No.1, 113-125.\bibitem{}
 C. Abdelkefi, Weighted function spaces and Dunkl transform. Mediterr.
J. Math. 9 (2012), 499-513 Springer Basel AG. \bibitem{} W.
Allegretto and P. O. Odiobala, Nonpositone elliptic problems in
$\mathbb{R}^{n}$, Proc. AMS 123 (1995) 533-541.\bibitem{} B. Amri
and M. Sifi, Riesz transforms For Dunkl transform, Annales
Math\'{e}matiques Blaise Pascal 19, 247-262 (2012).\bibitem{} J. S.
Bradley, Hardy inequalities with mixed norms, Canad. Math. Bull. 21
(1978), 405-408.
\bibitem{} A. P. Calder\'{o}n, Spaces between $L^1$ and $L^\infty$ and the theorem of Marcinkiewicz, Studia Math.
26 (1966), 273-299.\bibitem{} M. F. E. de Jeu, The Dunkl transform.
Inv. Math, 113 (1993), 147-162.
 \bibitem{} C. F. Dunkl, Differential-difference operators
associated to reflection groups. Trans. Amer. Math. Soc. 311, No1,
(1989), 167-183.
\bibitem{} C. F. Dunkl, Integral kernels with reflection group
invariance. Can. J. Math. 43, No 6, (1991), 1213-1227. \bibitem{} G.
H. Hardy, J. E. Littlewood and G. P\'{o}lya, Inequalities, 2nd ed.,
Cambridge Univ. Press, London/New York, 1952.\bibitem{} S. Hassani,
S. Mustapha, M. Sifi, Riesz potentials and fractional maximal
function for the Dunkl transform, J. Lie Theory 19 (4) (2009)
725-734.\bibitem{} T. Kurokawa and Y. Mizuta, On the order at
infinity of Riesz potentials, Hiroshima Math. J.9 (1979) 533-545.
\bibitem{} M. R\"osler and M. Voit, Markov processes with Dunkl
operators. Adv. in Appl. Math. 21, (1998), 575-643.\bibitem{} M.
R\"osler, Dunkl operators: theory and applications, in Orthogonal
polynomials and special functions (Leuven, 2002). Lect. Notes Math.
1817, Springer-Verlag (2003), 93-135.\bibitem{} M. R\"osler, A
positive radial product formula for the Dunkl kernel. Trans. Amer.
Math. Soc. 335, n°6, (2003), 2413-2438.\bibitem{} D. Siegel  and E.
Talvila, Pointwise growth estimates of the Riesz potential, Dynamics
of Continuous, Discrete and Impulsive System, 5 (1999), 185-194.
\bibitem{} S. Thangavelyu and Y. Xu, Convolution operator and
maximal function for Dunkl transform. J. Anal. Math. Vol. 97,
(2005), 25-56.\bibitem{} S. Thangavelu, Y. Xu, Riesz transform and
Riesz potentials for Dunkl transform, J. Comput. Appl. Math. 199
(2007) 181-195.\bibitem{} K. Trim\`eche, Paley-Wiener theorems for
the Dunkl transform and Dunkl translation operators. Integral
Transforms Spec. Funct. 13, (2002), 17-38.\bibitem{} J.F. van
Diejen, L. Vinet, Calogero-Sutherland-Moser Models. CRM Series in
Math. Phys., Springer-Verlag, 2000.
\end{thebibliography}
\end{document}